\documentclass[aop,citesort,MSNbibl,dvips]{arximspdf}

\doi{10.1214/09-AOP511}
\volume{38}
\issue{4}
\pubyear{2010}
\firstpage{1345}
\lastpage{1367}

\makeatletter

\newtheorem{prop}{Proposition}
\newtheorem{lem}{Lemma}
\newtheorem{theorem}{Theorem}
\newtheorem{thmref}{Theorem}
\newtheorem{lemref}{Lemma}

\newproclaim{Remarks}{Remarks}
\newproclaim{Remark}{Remark}
\newproclaim{Definition}{Definition}

\makeatother

\begin{document}
\begin{frontmatter}

\title{Uniform convergence of Vapnik--Chervonenkis classes
under ergodic sampling}
\runtitle{VC Classes under ergodic sampling}

\begin{aug}
\author[A]{\fnms{Terrence M.} \snm{Adams}\ead[label=e1]{tmadam2@tycho.ncsc.mil}} and
\author[B]{\fnms{Andrew B.} \snm{Nobel}\corref{}\thanksref{t1}\ead[label=e2]{nobel@email.unc.edu}}
\runauthor{T. M. Adams and A. B. Nobel}
\affiliation{Department of Defense and University of North Carolina at
Chapel Hill}
\address[A]{Department of Defense\\
9800 Savage Road\\
Suite 6513\\
Fort Meade, Maryland 20755\\
USA\\
\printead{e1}}
\address[B]{Department of Statistics\\
\quad and Operations Research\\
University of North Carolina\\
Chapel Hill, North Carolina 27599-3260\\
USA\\
\printead{e2}}
\end{aug}

\thankstext{t1}{Supported by NSF Grant DMS-04-06361.}

% HISTORY:
\received{\smonth{1} \syear{2008}}
\revised{\smonth{10} \syear{2009}}

% ABSTRACT
%
\begin{abstract}
We show that if $\mathcal{X}$ is a complete separable metric space and
$\mathcal{C}$ is a countable family of Borel subsets of $\mathcal{X}$
with finite VC dimension, then, for every stationary ergodic process
with values in $\mathcal{X}$, the relative frequencies of sets $C
\in\mathcal{C}$ converge uniformly to their limiting probabilities.
Beyond ergodicity, no assumptions are imposed on the sampling process,
and no regularity conditions are imposed on the elements of
$\mathcal{C}$. The result extends existing work of Vapnik and
Chervonenkis, among others, who have studied uniform convergence for
i.i.d. and strongly mixing processes. Our method of proof is new and
direct: it does not rely on symmetrization techniques, probability
inequalities or mixing conditions. The uniform convergence of relative
frequencies for VC-major and VC-graph classes of functions under
ergodic sampling is established as a corollary of the basic result for
sets.
\end{abstract}

% KEYWORDS
%

\end{frontmatter}

%s1 ###
\section{Introduction}\label{intro}

The strong law of large numbers and its extension to dependent
processes via the ergodic theorem
is one of the central results of probability theory. The strong law
connects sampling
and population-based quantities, and is one of the basic tools
for establishing the consistency of statistical inference procedures.
Uniform laws of large numbers extend the strong law by
guaranteeing the uniform convergence of averages to their limiting expectations
over a given family of functions.
Uniform laws of large numbers have been widely used and extensively studied
in a number of fields, including statistics, where they play a
foundational role
in the theory of empirical processes and machine learning. In the latter,
they underlie many results on consistency and rates of convergence
for classification and regression procedures.

The majority of the work on uniform laws of large numbers
to date has considered independent, identically distributed
samples, although there is also a substantial literature concerned with
dependent sequences satisfying
a variety of mixing conditions. The primary focus of this paper is the
uniform convergence of relative frequencies over a family of sets for general
ergodic processes.
In particular, we show that a sufficient condition for uniform
convergence in the
i.i.d. case, namely having finite Vapnik--Chervonenkis (VC) dimension,
is also sufficient to ensure uniform convergence in the ergodic case.
The VC dimension is a combinatorial quantity that describes
the ability of a collection of sets to pick apart finite subsets of points.
It can be defined without reference to metrics, epsilon-coverings,
metric entropies
or standard notions of vector space dimension.

Let $\mathbf{X}= X_1, X_2,\ldots$ be a stationary sequence of random variables
%defined on an underlying probability space $(\Omega, \F, \mathbb P)$
%and
taking values in a complete separable metric space $\mathcal{X}$ equipped
with its associated Borel sigma-field $\mathcal{S}$.
Under the standard definition,
$\mathbf{X}$ is ergodic if its invariant sigma-field is trivial
(cf. Definition 6.30 in Breiman \cite{Breim92}).
An equivalent, mixing-based definition of ergodicity can be formulated
as follows. For each $k \geq1$, let $\mathcal{S}^k$ denote the usual
product sigma-field on $\mathcal{X}^k$.
The process $\mathbf{X}$ is then ergodic if,
for each $k \geq1$ and every $A,B \in\mathcal{S}^k$,
%
%e1 ###
%
\begin{equation}
\label{ergmix}
\lim_{n \to\infty}
\frac{1}{n} \sum_{i=0}^{n-1} \mathbb P( X_1^k \in A, X_{i+1}^{i+k}
\in B )
\to
\mathbb P( X_1^k \in A ) \mathbb P( X_1^k \in B ) ,
\end{equation}
where $X_1^k$ denotes the $k$-tuple $(X_1,\ldots,X_k)$.
The condition simply states that, on
average, the present and the future of $\mathbf{X}$ become independent
as the
gap between them grows.

Suppose that $\mathbf{X}$ is ergodic. Here, and in what follows,
we let $X$ denote a
random variable independent of $\mathbf{X}$ and having the
same distribution as $X_1$.
For each set $C \in\mathcal{S}$,
the ergodic theorem ensures that the relative frequency
$m^{-1} \sum_{i=1}^m I_C(X_i)$ of $C$ converges almost surely
to the probability $\mathbb P(X \in C)$
as $m$ tends to infinity. Of interest here are families
of sets over which this convergence is uniform.
To this end, we define the random variables
%
%e2 ###
%
\begin{equation}
\label{gammamdef}
\Gamma_m(\mathcal{C}\dvtx\mathbf{X})
\stackrel{\triangle}{=}
\sup_{C \in\mathcal{C}}
\Biggl| \frac{1}{m} \sum_{i=1}^{m} I(X_i \in C) - \mathbb P(X \in C)
\Biggr|,\qquad m \geq1 .
\end{equation}
A countable
family $\mathcal{C}$ of Borel measurable sets is said to be a
\textit{Glivenko--Cantelli} class
for $\mathbf{X}$ if the relative frequencies of $C \in\mathcal{C}$
converge uniformly
to their limiting probabilities, in the sense that
%
%e3 ###
%
\begin{equation}
\label{UC}
\Gamma_m(\mathcal{C}\dvtx\mathbf{X})
\to0\qquad
\mbox{with probability one as } m \to\infty.
\end{equation}
Note that the uniformity here is over the family $\mathcal{C}$, not the
underlying sample
space; following standard usage, the term ``uniform convergence'' is
used rather
than the more traditional ``equiconvergence.''
The assumption that $\mathcal{C}$ is countable ensures that the supremum
in (\ref{UC}) is measurable. Uncountable families are discussed
briefly below.

Vapnik and Chervonenkis \cite{VC71} established necessary and
sufficient conditions
for (\ref{UC}) under i.i.d. sampling.
Their work provides a connection
between uniform convergence and the combinatorial complexity of a family
$\mathcal{C}$, where the latter is measured by the ability of the
family to break apart finite
sets of points. Let $\mathcal{C}$ be any collection of subsets of
$\mathcal{X}$ and let
$D \subseteq\mathcal{X}$ be any finite set of points.
The \textit{shatter coefficient} (or \textit{index}) of $\mathcal{C}$ with
respect to $D$ is defined by
%
%e4 ###
%
\begin{equation}
\label{shcdef}
S(D \dvtx\mathcal{C}) = | \{C \cap D \dvtx C \in\mathcal{C}\} |
\end{equation}
and is simply the number of distinct subsets of $D$ that
can be captured by sets $C \in\mathcal{C}$.
Clearly, $S(D \dvtx\mathcal{C}) \leq2^{|D|}$. When equality holds,
$\mathcal{C}$ is said to
shatter the set~$D$. The result of Vapnik and Chervonenkis can be stated
as follows.
\renewcommand{\thethmref}{A}
\begin{thmref}[(Vapnik and Chervonenkis \cite{VC71})]
\label{VC}
If $X_1, X_2,\ldots$ are i.i.d., then the uniform strong law
(\ref{UC}) holds if and only if
\[
\frac{1}{n} \log S(\{X_1, \ldots, X_{n} \} \dvtx\mathcal{C})
\rightarrow0
\]
in probability as $n$ tends to infinity.
\end{thmref}

In subsequent work, Vapnik and Chervonenkis \cite{VC81}
characterized uniform convergence
for classes of real-valued functions through the related notion of
metric entropy.
Talagrand \cite{Tala87} later provided a characterization of uniform
convergence in the
i.i.d. case
that strengthens these results and is focused on what happens when uniform
convergence fails. For nonatomic distributions, his results show that
(\ref{UC}) fails to hold if and only
if there is a set $A \in\mathcal{S}$ with $P(A) > 0$ such that,
for almost every realization of $\mathbf{X}$, the family
$\mathcal{C}$ shatters the set
$\{ X_{n_1}, X_{n_2}, \ldots\}$ consisting of those $X_i$ that lie in $A$.
\begin{Definition*}
The \textit{VC dimension} of a family $\mathcal{C}$, denoted here by
$\dim(\mathcal{C})$,
is the largest integer $k \geq1$ such that
$S(D \dvtx\mathcal{C}) = 2^k$ for some $k$-element subset $D$ of
$\mathcal{X}$.
If, for every $k \geq1$, the family $\mathcal{C}$ can shatter some
$k$-element set,
then $\dim(\mathcal{C}) = + \infty$.
\end{Definition*}

A family $\mathcal{C}$ is said to be a \textit{VC class} if
$\dim(\mathcal{C})$ is finite. The following combinatorial result of
Sauer provides
polynomial bounds on the shatter coefficients of VC classes in terms of
their combinatorial dimensions.
\renewcommand{\thelemref}{A}
\begin{lemref}[(Sauer \cite{Sau72})]
\label{Sauer}
If $\dim(\mathcal{C}) = V < \infty$, then
$S(D \dvtx\mathcal{C}) \leq\sum_{j=0}^{V} {m \choose j} \leq(m+1)^V$
for every $m \geq V$ and every
$D \subseteq\mathcal{X}$ of cardinality $m$.
\end{lemref}

It follows from Lemma \ref{Sauer} and Theorem \ref{VC} that if $V =
\dim(C) < \infty$,
then $\mathcal{C}$ is a Glivenko--Cantelli class for every i.i.d.
process $\mathbf{X}
$. Indeed,
one may establish an exponential inequality of the form
$\mathbb P( \Gamma_m(\mathcal{C}\dvtx\mathbf{X}) > t ) \leq c_1 (m+1)^V
e^{-c_2 m t^2}$ for
every $t > 0$ and $m \geq1$, where $c_1$ and $c_2$ are constants that are
independent of $m$, $\mathcal{C}$ and the distribution of $\mathbf
{X}$ (cf. \cite{DevLug01}).
The notions of VC class and VC dimension play a central role in modern
central limit and empirical process theory; see \cite
{VaaWel96,Dudley99} and
the references therein.

%s1.1 ###
\subsection{Principal result}\label{sec11}

In this paper, we show that the uniform strong law (\ref{UC}) holds
for VC classes under general ergodic sampling schemes. No mixing
conditions are imposed beyond ergodicity, and no conditions are imposed
on the elements of
$\mathcal{C}$. Under these circumstances, the convergence guaranteed
by the ergodic theorem can be arbitrarily slow and we cannot hope to obtain
distribution-free probability bounds like those discussed above for the
i.i.d. case.
Nevertheless, asymptotic results are still possible.
Our principal result is the following theorem; its proof can be found
in Sections
\ref{Special} and \ref{RPT} below.
\begin{theorem}
\label{ErgVC}
Let $\mathcal{X}$ be a complete separable metric space equipped with its
Borel measurable subsets $\mathcal{S}$ and
let $\mathcal{C}\subseteq\mathcal{S}$ be any countable family of
sets. If
$\dim(\mathcal{C}) < \infty$, then, for every stationary ergodic process
$\mathbf{X}= X_1,X_2, \ldots$ taking values in $(\mathcal
{X},\mathcal{S})$,
%
%e5 ###
%
\begin{equation}
\label{usl2}
\Gamma_m(\mathcal{C}\dvtx\mathbf{X}) =
\sup_{C \in\mathcal{C}}
\Biggl| \frac{1}{m} \sum_{i=1}^{m} I(X_i \in C) - \mathbb P(X \in C) \Biggr|
\to0 \qquad\mbox{w.p. 1},
\end{equation}
as $m$ tends to infinity. In other words, $\mathcal{C}$ is a Glivenko--Cantelli
class for every stationary ergodic process.
\end{theorem}

%s1.2 ###
\subsection{Uncountable families of sets}

The assumption that the family $\mathcal{C}$ is countable ensures
that the suprema $\Gamma_m(\mathcal{C}\dvtx\mathbf{X})$ are measurable and
is required for the construction of the isomorphism in Lemma
\ref{RedLem3}. In addition, countability of $\mathcal{C}$ is used in
the proof of
Proposition \ref{ErgVCweak} to ensure that no sample $X_i$ takes values
in the boundary of any set $C \in\mathcal{C}$.

Although it can be weakened in many cases
(see the discussion below), the assumption that $\mathcal{C}$ is countable
cannot be dropped altogether since it excludes somewhat pathological
examples that may arise in the dependent setting. To illustrate,
let $\mu$ be a nonatomic measure on $(\mathcal{X}, \mathcal{S})$ and
let $T\dvtx\mathcal{X}\to\mathcal{X}$ be an ergodic $\mu$-measure-preserving
bijection of $\mathcal{X}$. (More concretely, one may
take $T$ to be an irrational rotation of the unit circle
with its uniform measure.) Let $T^i$ denote the $i$-fold
composition of $T$ with itself if $i \geq1$, the $i$-fold
composition of $T^{-1}$ with itself if $i \leq-1$ and the
identity if $i = 0$.
For each $x \in\mathcal{X}$, let
$\mathcal{C}_{x} = \bigcup_{i=-\infty}^{\infty} \{ T^i x \}$
be the trajectory of $x$ under $T$ and
define the family $\mathcal{C}= \{ C_x \dvtx x \in\mathcal{X}\}$.
It is easy to see that for any two points $x_1, x_2 \in\mathcal{X}$, either
$C_{x_1} = C_{x_2}$ or
$C_{x_1} \cap C_{x_2} = \varnothing$, and so the VC dimension of
$\mathcal{C}$
equals one.
Now, let $X_i = T_i X_0$, where $X_0 \in\mathcal{X}$ is distributed
according to $\mu$. The process $\mathbf{X}= X_0, X_1, \ldots$
is then stationary and ergodic. Moreover, the $\mu$-measure of the
countable set $C_x$ is zero
for every $x$ and it is easy to see that
$\Gamma_m( \mathcal{C}\dvtx\mathbf{X}) = 1$ with probability one for
each $m \geq1$. Thus, (\ref{usl2}) fails to hold.

In spite of such negative examples,
Theorem \ref{ErgVC} can be extended in a straightforward way to
uncountable classes $\mathcal{C}$ under a natural approximation
condition. We
will call
an uncountable family $\mathcal{C}\subseteq\mathcal{S}$ ``nice''
for a given process
$\mathbf{X}$ if $\Gamma_m(\mathcal{C}\dvtx\mathbf{X})$ is measurable
for each $m \geq1$ and if,
for every $\varepsilon> 0$, there exists a countable subfamily
$\mathcal{C}_0 \subseteq\mathcal{C}$ such that
$\limsup_m \Gamma_m(\mathcal{C}\dvtx\mathbf{X}) \leq\limsup_m \Gamma
_m(\mathcal{C}_0\dvtx\mathbf{X}) +
\varepsilon$
with probability one. If $\mathcal{C}$ has finite VC dimension, then
(\ref{usl2}) holds for every ergodic process $\mathbf{X}$ such that
$\mathcal{C}$ is nice for $\mathbf{X}$.
%if there exists a countable sub-family $\C_0 \subseteq\C$ such that
%$\Gamma_m(\bX: \C) = \Gamma_m(\bX: \C_0)$ for each $m \geq1$ and
%every sample sequence of $\bX$.

Theorem \ref{ErgVC} can also be extended to the case in which the elements
of $\mathcal{C}$ belong to the completion of the Borel sigma-field of
$\mathcal{X}$
with respect to
the common distribution of the $X_i$.

%s1.3 ###
\subsection{Families of functions}

Theorem \ref{ErgVC} can be used to establish two related uniform
convergence results for families of functions. These results are
presented below.
In each case, the results can be extended to uncountable families
$\mathcal{F}$ under approximation conditions like those above for families
of sets.

A countable family $\mathcal{F}$ of Borel measurable functions
$f\dvtx\mathcal{X}\to\mathbb R$ is said to be a
\textit{Glivenko--Cantelli class} for
a stationary
ergodic process $\mathbf{X}$ if the relative frequencies
of functions in $f$ converge uniformly to their limiting expectations,
that is,
%
%e6 ###
%
\begin{equation}
\label{UC-fcn}
\Gamma_m(\mathcal{F}\dvtx\mathbf{X})
\stackrel{\triangle}{=}
\sup_{f \in\mathcal{F}}
\Biggl| \frac{1}{m} \sum_{i=1}^{m} f(X_i) - Ef(X) \Biggr|
\to0\qquad
\mbox{w.p. 1 as } m \to\infty.
\end{equation}
Here, we assume that the expectation $Ef(X)$ is well defined for
each $f \in\mathcal{F}$. Recall that a measurable function $F:
\mathcal{X}\to
[0,\infty)$
is said to be an \textit{envelope} for $\mathcal{F}$ if $|f(x)| \leq F(x)$
for each $x \in\mathcal{X}$ and $f \in\mathcal{F}$. In particular,
$\mathcal{F}$ is bounded
if it has constant envelope $F = M < \infty$.

%s1.3.1 ###
\subsubsection{VC-major classes}
\label{VCM}

Let $L_f(\alpha) = \{ x \dvtx f(x) \leq\alpha\}$ denote the $\alpha$
level set of a function
$f\dvtx\mathcal{X}\to\mathbb R$. A family of functions $\mathcal{F}$ is
said to be a \textit{VC-major class} if
\[
\dim_{\mathrm{VC}}(\mathcal{F}) = \sup_{\alpha\in\mathbb R} \dim\bigl( \{
L_f(\alpha) \dvtx f
\in\mathcal{F}\} \bigr)
\]
is finite. The following result is established in Section \ref{VCMajorProof}.
\begin{prop}
\label{VCMajor}
Let $\mathcal{F}$ be a countable family of Borel measurable functions
$f\dvtx\mathcal{X}
\to\mathbb R$ with
envelope $F$. If $\mathcal{F}$ is a VC-major class, then (\ref
{UC-fcn}) holds
for every
stationary ergodic process $\mathbf{X}$ such that $EF(X)$ is finite.
\end{prop}

%s1.3.2 ###
\subsubsection{VC-graph classes}
\label{VCG}

The graph of a function $f\dvtx\mathcal{X}\to\mathbb R$ is the set $G_f
\subseteq\mathcal{X}
\times\mathbb R$
defined by $G_f = \{ (x,s) \dvtx
%x \in\X\mbox{ and }
0 \leq s \leq f(x) \mbox{ or } f(x) \leq s \leq0 \}$.
A family $\mathcal{F}$ of functions $f\dvtx\mathcal{X}\to\mathbb R$ is
said to be a \textit{VC-graph
class} (Pollard \cite{Poll84}) if
\[
\dim_G(\mathcal{F}) = \dim(\{ G_f \dvtx f \in\mathcal{F}\})
\]
is finite. The
following result is established in Section \ref{VCGraphProof}.
\begin{prop}
\label{VCGraph}
Let $\mathcal{F}$ be a bounded, countable family of Borel measurable functions
$f\dvtx\mathcal{X}\to\mathbb R$. If $\mathcal{F}$ is a VC-graph class,
then (\ref{UC-fcn})
holds for every
stationary ergodic process $\mathbf{X}$.
\end{prop}

%s1.4 ###
\subsection{Related work}

Steele \cite{Stee78} used subadditive ergodic theory to establish
that both $\Gamma_m(\mathcal{C}\dvtx\mathbf{X})$ [see (\ref
{gammamdef})] and the entropy
$n^{-1} \log S(\{X_1,\ldots,X_n \}\dvtx\mathcal{C})$ [see (\ref{shcdef})]
converge with probability
one to nonnegative constants whenever $\mathbf{X}$ is ergodic.
In addition, he obtained refined
necessary and sufficient conditions for uniform strong laws in the
i.i.d. case.
Nobel \cite{Nob95} showed that the conditions of Theorem \ref{VC} and
Talagrand \cite{Tala87} do not characterize uniform convergence for ergodic
processes and, in particular, that standard random entropy conditions
do not ensure uniform convergence in the ergodic case.

Yukich \cite{Yuk86} established rates of convergence for $\Gamma
_m(\mathcal{F}
\dvtx\mathbf{X})$
when $\mathbf{X}$ is $\phi$-mixing
and $\mathcal{F}$ satisfies suitable bracketing entropy conditions.
Yu \cite{Yu94} extends these results to
$\beta$-mixing (absolutely regular) processes $\mathbf{X}$ and
classes $\mathcal{F}$ satisfying metric entropy conditions.
(See Bradley \cite{Brad86} for more on $\phi$- and $\beta$-mixing
conditions.)
For VC classes $\mathcal{C}$, the results of Yu imply the uniform law
(\ref
{usl2}) when
the mixing coefficients $\beta_k$
decrease as $k^{-r}$ for some $r > 0$.
Work of Pe{\v{s}}kir and Yukich \cite{PesYuk94} extends this conclusion to
$\beta$-mixing
processes with $\beta_k = (\log k)^{-2}$.
%They extended these results to families of processes, measure
%preserving flows, and sequences generated by the iteration of a
%nonlinear operator.

Nobel and Dembo \cite{NobDem93} showed that one may extend uniform
strong laws
from i.i.d. processes to $\beta$-mixing processes with the same
one-dimensional marginal distribution. Their result implies that (\ref{usl2})
holds for any VC class $\mathcal{C}$ and any $\beta$-mixing process
$\mathbf{X}$.
Peligrad \cite{Peli01} established an analogous result for processes satisfying
a modified $\phi$-mixing condition.
Karandikar and Vidyasagar \cite{KarVid02} extended the results of
\cite{NobDem93}
to families of processes and established rates of
convergence depending on the behavior of the mixing coefficients.

Extending earlier work of Hoffmann-J{\o}rgensen \cite{Hoff84}
in the i.i.d. case, Pe{\v{s}}kir and Weber \cite{PesWeb93}
show that the uniform ergodic
theorem (\ref{UC-fcn}) holds if and only if the family $\mathcal{F}$
is, in
their terminology,
eventually totally bounded in mean. They also note
the equivalence of different notions of convergence, as in Steele's
work.
Pe{\v{s}}kir \cite{Pes98} investigated conditions for uniform mean square ergodic
theorems for families of weak-sense stationary processes.
%under conditions on their spectral representations.

Andrews \cite{And87} investigated sufficient conditions under which laws
of large numbers can be extended from individual functions to
classes of functions, with particular emphasis on stochastically
equicontinuous classes indexed by
totally bounded parameter spaces. The bibliography of
his paper provides a good overview of related work.

%The next section presents several technical preliminaries.
%The proof of Theorem \ref{ErgVC}
%is contained in Section \ref{Proof}. Section \ref{PRCF} contains the
%proofs
%of Propositions \ref{VCMajor} and \ref{VCGraph}.

%s1.5 ###
\subsection{Overview}

In the absence of independence or standard uniform mixing conditions, a
direct approach to Theorem \ref{ErgVC} using symmetrization and
exponential-type inequalities, or a more indirect approach carried out by
coupling with the independent case, does not
appear to be possible. Instead, we establish, without reference to
independence or mixing conditions, the contrapositive of Theorem~\ref{ErgVC}:
if the relative frequencies
of sets $C \in\mathcal{C}$ fail to converge uniformly, then, for each
$L \geq
1$, we
can find $L$ points $x_1,\ldots,x_L \in\mathcal{X}$ that are
shattered by $\mathcal{C}$
and, consequently, $\dim(\mathcal{C}) = \infty$.
For this, we require only the almost sure convergence guaranteed by
the ergodic theorem for individual sets.
Rather than working directly with the shatter coefficients $S(\cdot\dvtx
\mathcal{C})$,
we consider joins (partitions) generated by finite subcollections of
$\mathcal{C}$,
which are defined in Section \ref{Special} below.

In the next section, we begin with a special case of Theorem \ref{ErgVC}
in which $\mathcal{X}= [0,1]$, each $X_i$ is uniformly
distributed on $\mathcal{X}$ and each element of $\mathcal{C}$ is
equal to a finite
union of intervals.
This preliminary result, which is the core of the paper,
is contained in Proposition \ref{ErgVCweak}.
The general case of Theorem \ref{ErgVC} is established in
Section \ref{RPT} using Proposition \ref{ErgVCweak} and a series of
three reductions.
The first reduction (contained in Lemma \ref{RedLem1})
shows that it is enough to consider
processes $\mathbf{X}$ for which the marginal distribution of the
$X_i$ is nonatomic.
The second reduction maps the random variables in $\mathbf{X}$ and
the elements of $\mathcal{C}$ to the unit
interval with Lebesgue measure via a standard measure space isomorphism.
The final reduction (contained in Lemma \ref{RedLem3})
makes use of an additional measure space
isomorphism that maps each element of $\mathcal{C}$ into a set that is
equal, up to a set of measure zero, to a finite union of intervals.

%s2 ###
\section{Classes containing finite unions of intervals}
\label{Special}

In this section, we establish a version of Theorem \ref{ErgVC} in which
$\mathcal{X}= [0,1]$ and each element of $\mathcal{C}$ is a finite
union of intervals.
In the proof, we work with the joins of selected members of $\mathcal
{C}$, which
act as surrogates for the more commonly used shatter coefficients.
\begin{Definition*} The \textit{join} of $k$ sets $A_1,\ldots, A_k
\subseteq\mathcal{X}$,
denoted $J = \bigvee_{i=1}^k A_i$, is the collection of all
\textit{nonempty} intersections $\tilde{A}_1 \cap\cdots\cap\tilde{A}_k$,
where $\tilde{A}_i \in\{ A_i, A_i^c \}$ for $i = 1,\ldots,k$. Note
that $J$
is a partition of $\mathcal{X}$.
The join of $A_1,\ldots, A_k$ is said to be \textit{full} if it has
(maximal) cardinality $2^k$.
%It is a fat join if the measure of each cell in the join is positive.
%In the proof of
%Theorem \ref{ErgVC}, we will use the fat join notion.
%Note that any fat join is a full join, and arbitrarily large
%full joins produce arbitrarily large VC-dimension.
% \textit{i.e.}, $|\bigvee_{i=1}^k A_i| = 2^k$.
\end{Definition*}

The next lemma makes an elementary connection between full joins
and the VC dimension.
% ABN: Add Citation
A similar result appears in \cite{Mat02} as Lemma 10.3.4.
We include a short proof here for completeness.
\begin{lem}
\label{join}
Let $\mathcal{C}$ be any collection of subsets of $\mathcal{X}$. If,
for some $k \geq1$,
there exists a collection $\mathcal{C}_0 \subseteq\mathcal{C}$ of
$2^k$ sets having
a full join, then VC-$\dim(\mathcal{C}) \geq k$.
\end{lem}
\begin{pf}
Indexing the elements of $\mathcal{C}_0$ in an arbitrary
manner by
subsets of $[k] : = \{1,\ldots,k\}$, we may write
$\mathcal{C}_0 = \{ C(U) \dvtx U \subseteq[k] \}$. For $i = 1,\ldots
,k$, let
$x_i$ be
any element of the intersection
\[
\biggl( \bigcap_{ U \subseteq[k], i \in U } C(U) \biggr)
\cap
\biggl( \bigcap_{ U \subseteq[k], i \notin U } C(U)^c \biggr) ,
\]
which is nonempty by assumption. For each subset
$V \subseteq[k]$, it is easy to see that $x_i \in C(V)$ if and only
if $i \in V$. Thus, $\mathcal{C}_0$, and hence $\mathcal{C}$, shatters
$\{x_1,\ldots,x_k\}$.

Now, let $\mathbf{X}= X_1, X_2, \ldots$ be a stationary ergodic process
defined on a probability space $(\Omega, \mathcal{F}, \mathbb P)$,
such that
each $X_i$ takes values in $[0,1]$, is Borel measurable and
has distribution equal to the Lebesgue measure $\lambda(\cdot)$
on $[0,1]$.
\end{pf}
\begin{prop}
\label{ErgVCweak}
Let $\mathcal{C}_0$ be a countable family of subsets of $[0,1]$,
each of whose elements is a finite union of intervals. Suppose that
\[
\limsup_{m} \Gamma_m(\mathcal{C}_0\dvtx\mathbf{X}) > 0
\]
with positive probability. Then, for each integer $L \geq1$, there
exist sets $D_1,\break D_2, \ldots, D_L \in\mathcal{C}_0$ such that the join
$K_L = D_1 \vee D_2 \vee\cdots\vee D_L$ is full and each cell
of $K_L$ has positive Lebesgue measure.
\end{prop}
\begin{Remarks*} It follows from Lemma \ref{join} that the family
$\mathcal{C}_0$ in Proposition \ref{ErgVCweak} has infinite VC dimension.
The additional fact that each cell of the joins has positive measure
will be
needed in the proof of Theorem \ref{ErgVC} as we may then ignore
sets of measure zero that arise in the application of Lemma
\ref{RedLem3}. The assumption that $C \in\mathcal{C}_0$ is a finite
union of intervals guarantees that its boundary has
Lebesgue measure zero. Excluding such boundary points
from the process $\mathbf{X}$ plays an important role in the final
part of the proof of Lemma \ref{fulljoin}.
\end{Remarks*}
\begin{pf*}{Proof of Proposition \protect\ref{ErgVCweak}}
In what follows, we will need to examine the difference
between the relative frequency and probability of
subsets of the unit interval. To this end, for each
$\omega\in\Omega$, each $A \subseteq[0,1]$ and
each $m \geq1$, we define
%
%e7 ###
%
\begin{equation}
\label{discrep}
\Delta^\omega(A\dvtx m)
\stackrel{\triangle}{=}
\Biggl| \frac{1}{m} \sum_{i=1}^{m} I \bigl( X_i(\omega) \in A \bigr) - \lambda(A ) \Biggr|
\end{equation}
to be the discrepancy of $A$ with respect to the first $m$ elements
of the sample sequence $X_i(\omega)$.
Let $B^o$, $\overline{B}$ and $\partial B = \overline{B} \setminus
B^o$ denote,
respectively, the interior, closure and boundary of a set $B \subseteq[0,1]$.

For $n \geq1$, let $\mathcal{D}_n = \{ [ k 2^{-n}, (k+1) 2^{-n}] \dvtx0
\leq k
\leq2^n-1 \}$
be the set of closed dyadic intervals of order $n$. Let
$\mathcal{D}$ be the union of the families $\mathcal{D}_n$ and let
$\mathcal{C}= \mathcal{C}_0 \cup\mathcal{D}$.
Then $\mathcal{C}$ and the set $A_0 = \bigcup_{C \in\mathcal{C}}
\partial C$
of all endpoints of elements of $\mathcal{C}$ are countable. In particular,
$\lambda(A_0) = 0$. By removing a $\mathbb P$-null set of
outcomes from our sample space, we can, and do, assume that
$X_i(\omega) \in A_0^c$ for every $i \geq1$ and every
$\omega\in\Omega$.

Recalling the definition (\ref{gammamdef}), we see that
$\Gamma_m(\mathcal{C}\dvtx\mathbf{X}) \geq\Gamma_m(\mathcal{C}_0 \dvtx
\mathbf{X})$
and, therefore, $\limsup_m \Gamma_m(\mathcal{C}\dvtx\mathbf{X}) > 0$
with positive probability.
In particular, there exists an \mbox{$\eta> 0$} and a set $E \in\mathcal{F}$
with $\mathbb P(E)>0$ such that
%
%e8 ###
%
\begin{equation}
\label{etadef}
\limsup_{m \to\infty}
\Bigl[
\sup_{C \in\mathcal{C}}
\Delta^\omega(C \dvtx m)
% | \frac{1}{m} \sum_{i=1}^{m} I( X_i(\omega) \in C ) - \lambda(C) |
\Bigr]
> \eta\qquad
\mbox{for each $\omega\in E$.}
\end{equation}
(Using the results of Steele \cite{Stee78} or, alternatively,
the invariance of $E$, it follows that that $\mathbb P(E) = 1$,
but we do not require this stronger result here.)
Fix $0 < \delta\leq\min\{ \eta/12, \mathbb P(E) \}$.

%ABN: NEW PARA
The remainder of the proof proceeds as follows. We first construct
a sequence of ``splitting sets'' $R_1, R_2, \ldots\subseteq[0,1]$, in
stages, from the sets in $\mathcal{C}$.
At the $k$th stage, the splitting set $R_k$ is obtained from
a sequential procedure that makes use of the splitting sets
$R_1,\ldots, R_{k-1}$ produced at previous stages. Once obtained,
the splitting sets are used to identify, for any $L \geq1$, a collection
of $L$ sets in $\mathcal{C}$ that have full join and
it is easy to show that at most one
member of such a collection can come from $\mathcal{D}$.
The final step of the proof requires that we keep track of the
process by which each splitting set $R_k$ is produced; this
requirement is reflected in the notation adopted below.

\subsection*{Construction of $R_1$}
We first choose a sequence of sets
$C_1, C_2, \ldots\in\mathcal{C}$ in such a way that a significant
fraction of
the cells in the join of $C_1,\ldots, C_n$ will intersect both
$C_{n+1}$ and its complement.
Let $C_1$ be any set in $\mathcal{C}$. Suppose
that $C_1,\ldots,C_n \in\mathcal{C}$ have already been selected and
we wish
to choose $C_{n+1}$. Let $J_n = \mathcal{D}_n \vee C_1 \vee\cdots
\vee C_n$ be
the join of the previously selected sets and the dyadic intervals of
order $n$.
Since the process $\mathbf{X}$ is ergodic and $J_n$ is finite,
there exists an integer $M$ and a set $F$ with $\mathbb P(F) > 1-\delta$
such that
%
%e9 ###
%
\begin{equation}
\label{deltacond}
\Delta^\omega(A \dvtx m) \leq\delta\lambda(A)\qquad
\mbox{for each $\omega\in F$, $A \in J_n$ and $m \geq M$}.
\end{equation}
As $\delta< \mathbb P(E)$, the set $E \cap F$ has positive $\mathbb P$-measure
and is therefore nonempty. Let $\omega_{n+1}$ be any point in $E \cap F$.
As $\omega_{n+1} \in E$,
it follows from (\ref{etadef}) that there exists a set $C_{n+1} \in
\mathcal{C}
$ and an
integer $m_{n+1} \geq M$ such that
$\Delta^{\omega_{n+1}} (C_{n+1} \dvtx m_{n+1}) > \eta$.
From $C_{n+1}$, one may construct the join
$J_{n+1} = \mathcal{D}_{n+1} \vee C_1 \vee\cdots\vee C_{n+1}$ and then
select $C_{n+2}$ in the same manner as $C_{n+1}$. Continuing in this
fashion, we obtain joins $J_{n+1}, J_{n+2}, \ldots$
and sets $C_{n+2}, C_{n+3}, \ldots\in\mathcal{C}$. We note that the sample
points $\omega_n$ may vary from step to step and that there is
no requirement that $m_{n+1}$ be greater than $m_n$.

The choice of the set $C_{n+1}$ ensures that it cannot be well
approximated by a
union of elements of $J_n$ or, equivalently, that the collection of
cells $A \in J_n$ containing points in $C_{n+1}$ and $C_{n+1}^c$
must have nonvanishing probability.
%This latter property lies at the root of the proof.
To make this idea precise, we define the family
\[
H_{n}
=
\biggl\{ A \in J_n\dvtx
\Delta^{\omega_{n+1}} (A \cap C_{n+1} \dvtx m_{n+1})
> \frac{\eta}{2} \lambda(A) \biggr\}.
\]
%
%where $\omega\in E \cap F$ is the point used to define $C_{n+1}$.
The next lemma shows that
the elements of $H_{n} \subseteq J_n$ occupy a nonvanishing
fraction of the unit interval.
\begin{lem}
\label{Glem1}
If $G_{n} = \bigcup H_{n}$ is the union of the sets $A \in H_{n}$, then
$\lambda(G_{n}) \geq\eta/ 6$.
\end{lem}
\begin{pf}
Let $\omega= \omega_{n+1}$, $C = C_{n+1}$ and $m = m_{n+1}$.
By decomposing\break $\Delta^\omega(C\dvtx m)$ among the
elements of $J_n$, we obtain the following bound:
%
%e10 ###
%
\begin{eqnarray}\label{rhseta}
\eta
& < &
\Delta^\omega(C \dvtx m)
\leq
\sum_{A \in J_n} \Delta^\omega(C \cap A \dvtx m)
\nonumber\\[-8pt]\\[-8pt]
& = &
\sum_{A \in H_{n}} \Delta^\omega(C \cap A \dvtx m)
+
\sum_{A \in J_n \setminus H_{n}} \Delta^\omega(C \cap A \dvtx
m).\nonumber
\end{eqnarray}
By definition of $H_{n}$, the second term in (\ref{rhseta}) is at most
\[
\sum_{A \in J_n \setminus H_{n}} \frac{\eta}{2} \lambda(A)
% = \frac{\eta}{2} \lambda(\cup(J_n \setminus H_n))
\leq
\frac{\eta}{2} .
\]
Moreover, the first term in (\ref{rhseta}) can be bounded as follows:
\begin{eqnarray*}
&&
\sum_{A \in H_{n}} \Delta^\omega(C \cap A \dvtx m) \\
&&\qquad \leq
\sum_{A \in H_{n}}
\frac{1}{m_{n+1}} \sum_{i=1}^{m} I\bigl( X_i(\omega) \in C \cap A \bigr)
+
\sum_{A \in H_{n}} \lambda(C \cap A) \\
&&\qquad \leq
\sum_{A \in H_{n}}
\frac{1}{m} \sum_{i=1}^{m} I\bigl( X_i(\omega) \in C \cap A \bigr)
+
\lambda(G_{n}) \\
&&\qquad \leq
\sum_{A \in H_{n}}
\Delta^\omega(A \dvtx m)
+
2 \lambda(G_{n}) \\
&&\qquad \leq
(\delta+ 2) \lambda(G_{n})
\leq
3 \lambda(G_{n}) ,
\end{eqnarray*}
where the penultimate inequality follows from (\ref{deltacond}) and the
fact that $\omega_{n+1} \in F$.
Combining the final expressions in the three preceding displays yields
the result.
%\rightqed
\end{pf}

Let the sets $G_n = \bigcup H_{n}$, $n \geq1$, be derived from the
inductive procedure described above. For each $n \geq1$, define a
sub-probability
measure $\lambda_n (B) = \lambda(B \cap G_n)$ on $([0,1],\mathcal
{B})$. The
collection $\{ \lambda_n \}$ is necessarily tight and therefore has a
subsequence
$\{ \lambda_{n_r} \}$ that converges weakly to a sub-probability $\nu
_1$ on
$([0,1],\mathcal{B})$, in the sense that
$\int_0^1 g \,d\lambda_{n_r} \to\int_0^1 g \,d\nu_1$ as $r \to\infty$
for every (bounded) continuous function $g\dvtx[0,1] \to\mathbb R$.
It is easy to see that $\nu_1$ is absolutely continuous with respect
to $\lambda$ and that
\[
\nu_1([0,1])
\geq
\limsup_{r \to\infty} \lambda_{n_r} ([0,1])
\geq
\eta/ 6 .
\]
In particular, the Radon--Nikodym derivative $d\nu_1 / d\lambda$ is
well defined
and is bounded above by 1.
Define $R_1 = \{ x\dvtx(d\nu_1 / d\lambda)(x) > \delta\}$. From the previous
remarks, it follows that
%
%e11 ###
%
\begin{eqnarray}\label{Rlbd}
\frac{\eta}{6}
& \leq&
\nu_1([0,1])
=
\int_0^1 \frac{d\nu_1}{d\lambda} \,d\lambda
=
\int_{R_1} \frac{d\nu_1}{d\lambda} \,d\lambda
+ \int_{R_1^c} \frac{d\nu_1}{d\lambda} \,d\lambda\nonumber\\[-8pt]\\[-8pt]
& \leq&
\int_{R_1} 1 \,d\lambda+ \int_{R_1^c} \delta \,d\lambda
\leq
\lambda(R_1) + \delta.\nonumber
\end{eqnarray}
As $\delta< \eta/12$ by assumption, we conclude that
$\lambda(R_1) \geq\eta/12 > 0$.

\subsection*{Construction of $R_k$ for $k \geq2$}
The splitting sets $R_2, R_3, \ldots$ are defined in order, following the
general iterative procedure used to construct $R_1$.
The critical difference between the first and subsequent stages is that
the sets
$R_1, \ldots, R_{k-1}$ produced at stages $1$ through $k-1$ are included
in the join used at stage $k$ to define~$R_{k}$.
In what follows, let $C_k(n)$, $J_k(n)$,
$\omega_k(n)$, $m_k(n)$, $H_k(n)$ and $G_k(n)$ denote
the quantities appearing at the $n$th step of the $k$th stage.
%, and let $n_k(l)$ denote the subsequence used to
In particular, let $C_1(n) = C_n$, $n \geq1$, be the elements
of $\mathcal{C}$ considered in stage 1 and define $J_1(n)$, $\omega_1(n)$,
$m_1(n)$, $H_1(n)$ and $G_1(n)$ in a similar fashion.

Suppose that stages $1$ through $k-1$ have been completed
and that we wish to construct the splitting set $R_k$ at stage $k$.
Let $C_k(1)$ be any element of $\mathcal{C}$
and suppose that $C_k(2), \ldots, C_k(n)$ have already been
selected. We define the join
%
%e12 ###
%
\begin{equation}
\label{joindef}
J_k(n) = \mathcal{D}_n \vee\bigvee_{j=1}^{k-1} R_j \vee\bigvee_{i=1}^n
C_k(i).
\end{equation}
By the ergodic theorem, there exists an integer $M$ and a set $F$
with $\mathbb P(F) > 1-\delta$ such that (\ref{deltacond}) holds with $J_n$
replaced by $J_k(n)$.
As before, it follows from these inequalities and
(\ref{etadef}) that there exists a
sample point $\omega_k(n+1) \in E \cap F$, a
set $C_k(n+1) \in\mathcal{C}$ and an integer $m_k(n+1) \geq M$
such that
%
%e13 ###
%
\begin{equation}
\label{deltaa}
\Delta^{\omega_k(n+1)} \bigl( A \dvtx m_k(n+1) \bigr)
\leq
\delta\lambda(A)\qquad
\mbox{for each } A \in J_k(n),
\end{equation}
and, simultaneously,
%
%e14 ###
%
\begin{equation}
\label{deltaeta}
\Delta^{\omega_k(n+1)} \bigl( C_k(n+1) \dvtx m_k(n+1) \bigr) > \eta.
\end{equation}
Using these quantities, we define the family
%
%e15 ###
%
\begin{equation}
\label{hkdef}\qquad
H_k(n)
=
\biggl\{ A \in J_k(n) \dvtx
\Delta^{\omega_k(n+1)}\bigl( C_k(n+1) \cap A \dvtx m_k(n+1) \bigr)
> \frac{\eta}{2} \lambda(A) \biggr\}
\end{equation}
and let $G_k(n) = \bigcup H_k(n)$ be the union of the elements of $H_k(n)$.

Defining $J_k(n+1)$ as in (\ref{joindef}) and continuing in the
same fashion,
we obtain a sequence $C_k(n+2), C_k(n+3),\ldots\in\mathcal{C}$ and
a corresponding sequence of sets
$G_k(n+1), G_k(n+2),\ldots\subseteq[0,1]$.
Lemma \ref{Glem1} ensures that $\lambda(G_k(n)) \geq\eta/ 6$ for
each $n \geq1$. As before, there is a sequence of integers
$n_k(1) < n_k(2) < \cdots$ such that the measures
$\lambda(B \cap G_k(n_k(r)))$ converge weakly as $r \to\infty$
to a sub-probability measure $\nu_k$ on $([0,1], \mathcal{B})$ that
is absolutely
continuous with respect to $\lambda(\cdot)$.
Define $R_k = \{ x\dvtx(d\nu_k / d\lambda)(x) > \delta\}$. The argument
in (\ref{Rlbd}) shows that $\lambda(R_k) \geq\eta/12$.
The arguments below require that we consider density points
of the splitting sets. With this in mind, for $k \geq1$, let
\[
\tilde{R}_k
=
\biggl\{ x \in R_k \dvtx\lim_{\alpha\to0}
\frac{ \lambda((x-\alpha,x+\alpha) \cap R_k)}
{2 \alpha} = 1 \biggr\}
\]
be the set of Lebesgue points of $R_k$.
By standard results on differentiation of integrals (cf. Theorem 31.3 of
Billingsley \cite{Bill95}),
$\lambda(\tilde{R}_k) = \lambda(R_k) \geq\eta/ 12$.
The sets $\tilde{R}_k$ are used to construct full joins in the next step
of the proof.

\subsection*{Construction of full joins}
Fix an integer $L \geq2$. As the measures of the sets $\tilde{R}_k$
are bounded away from
zero, there exist positive integers $k_1 < k_2 < \cdots< k_{L}$ such that
$\lambda(\bigcap_{j=1}^{L} \tilde{R}_{k_j}) > 0$. Define the intersections
\[
Q_r = \bigcap_{j=1}^{L-r} \tilde{R}_{k_j}
\]
for $r = 0, 1, \ldots, L-1$.
Note that $Q_0 \subseteq Q_1 \subseteq\cdots\subseteq Q_{L-1}$.
Recall that $B^o$, $\overline{B}$ and $\partial B$ denote, respectively,
the interior, closure and boundary of a set $B \subseteq[0,1]$.
\begin{lem}
\label{fulljoin}
There exist sets $D_1, D_2, \ldots, D_{L-1} \in\mathcal{C}$ such
that for each
$l = 1,\ldots,L-1$, the join
$K_l = D_1 \vee D_2 \vee\cdots\vee D_l$ satisfies
$|K_l| = 2^l$ and for each $B \in K_l$, the intersection
$B^o \cap Q_l$ is nonempty. In particular, each cell
of $K_l$ has positive Lebesgue measure.
%has positive $\lambda$-measure.
\end{lem}
\begin{pf}
We establish the result by induction on $l$, beginning
with the case $l = 1$. In particular, we show that there exists a
set $D_1 \in\mathcal{C}$ such that $D_1^o \cap Q_1$ and
$(D_1^c)^o \cap Q_1$ are nonempty. To this end,
we choose $x_1 \in Q_0$, which is nonempty by assumption,
and let $\varepsilon= \delta/ 2 (\delta+ 1)$. By
definition of the sets $\tilde{R}_{k_j}$, there exists
$\alpha_1 > 0$ such that the interval
$I_1 \stackrel{\triangle}{=}(x_1 - \alpha_1, x_1 + \alpha_1)$ satisfies
%
%e16 ###
%
\begin{equation}
\label{Int1prop}
\lambda(I_1 \cap Q_0)
\geq
(1 - \varepsilon) \lambda(I_1)
=
2 \alpha_1 (1 - \varepsilon) .
\end{equation}
To simplify notation, let $\kappa= k_L$. It follows from
the last display and the definition of $R_\kappa\supseteq Q_0$ that
%
%e17 ###
%
\begin{equation}
\label{nukineq}
\nu_{\kappa}(I_1 \cap R_{\kappa})
=
\int_{I_1 \cap R_\kappa} \frac{ d\nu_{\kappa} }{d\lambda}\,
d\lambda
>
\delta\lambda(I_1 \cap R_\kappa)
\geq
2 \alpha_1 (1 - \varepsilon) \delta.
\end{equation}
Now, let $\{ n_\kappa(r) \dvtx r \geq1 \}$ be the subsequence used to define
the sub-probability $\nu_\kappa$.
As $I_1$ is an open set, the portmanteau theorem and (\ref{nukineq})
imply that
\[
\liminf_{r \to\infty}
\lambda\bigl( I_1 \cap G_{\kappa}(n_\kappa(r)) \bigr)
\geq
\nu_{\kappa} (I_1)
\geq
\nu_{\kappa} (I_1 \cap R_\kappa)
>
2 \alpha_1 (1 - \varepsilon) \delta.
\]
Choose $r$ sufficiently large so that
$\lambda( I_1 \cap G_{\kappa}(n_\kappa(r)) ) >
2 \alpha_1 (1 - \varepsilon) \delta$ and
$2^{ - n_\kappa(r) } < \delta\alpha_1 / 4$. We require
the following subsidiary lemma.
\end{pf}
\begin{lem}
\label{Alem}
There exists a set $A \in H_\kappa(n_\kappa(r))$ such that
$A \subseteq I_1$ and $\lambda(A \cap Q_1) > 0$. Moreover,
$A$ is contained in $Q_1$.
\end{lem}
\begin{pf}
Let $G = G_{\kappa}(n_\kappa(r))$.
The choice of $n_\kappa(r)$ ensures that
\begin{eqnarray*}
(1 - \varepsilon) \delta\lambda(I_1)
& < &
\lambda( I_1 \cap G ) \\
& = &
\lambda( I_1 \cap Q_1 \cap G )
+
\lambda( I_1 \cap Q_1^c \cap G ) \\
& \leq&
\lambda( I_1 \cap Q_1 \cap G )
+
\lambda( I_1 \cap Q_1^c ) \\
& \leq&
\lambda( I_1 \cap Q_1 \cap G )
+
\varepsilon\lambda(I_1),
\end{eqnarray*}
where the first inequality follows from our choice of $r$
and the final inequality follows from (\ref{Int1prop}) together with the
fact that $Q_0 \subseteq Q_1$.
The last display and the definition of $\varepsilon$ imply
that $\lambda( I_1 \cap Q_1 \cap G ) \geq\delta\alpha_1$.
As the collection of sets used to define $J_\kappa(n_\kappa(r))$
includes the dyadic intervals of order $n_\kappa(r)$, each element
$A$ of the join has diameter (and Lebesgue measure) bounded by
$2^{ - n_\kappa(r) } < \delta\alpha_1 / 4$. These last two
inequalities imply that
\[
% =
\delta\alpha_1
\leq
\lambda( I_1 \cap Q_1 \cap G )
\leq
\sum_{A} \lambda(Q_1 \cap A) + 2 \frac{\delta\alpha_1}{4} ,
\]
where the sum is over sets $A \in H_\kappa(n_\kappa(r))$ such that
$A \subseteq I_1$. In particular, it is clear that the sum is necessarily
positive and the first part of the claim follows. Note that
$A \in H_\kappa(n_\kappa(r))$ implies that $A \in J_\kappa(n_\kappa(r))$.
Thus, the inclusion of the sets $R_1,\ldots, R_{\kappa- 1}$ in the
join ensures that
$A$ is contained in either $R_{k_j}$ or $R_{k_j}^c$, but not both,
for each $j=1,\ldots,L-1$.
If $\lambda(A \cap Q_1) > 0$, then, necessarily,
$A \cap Q_1 \neq\varnothing$, and the containment relations
imply that $A \subseteq Q_1$. This completes the proof of
Lemma \ref{Alem}.
%Also note that $\lambda(A) > 0$.
\end{pf}

Let $D_1 = C_\kappa(n_\kappa(r)+1) \in\mathcal{C}$, where $r$ is the
index appearing in Lemma \ref{Alem}. Recall that $D_1$ is a
finite union of intervals and that no random variables
$X_i$ take values in the finite set $\partial D_1$. In addition,
$\partial D_1$ has
Lebesgue measure zero.
Let $A \in H_\kappa(n_\kappa(r))$ be the set identified in
Lemma \ref{Alem} and note that $\lambda(A) > 0$.
We argue by contradiction that $A$ (and therefore $Q_1$)
has nonempty intersection
with the interiors of $D_1$ and~$D_1^c$. Suppose, first, that
$A \cap D_1^o = \varnothing$. In this case,
\[
\Delta^\omega(A \cap D_1 \dvtx m)
=
\Delta^\omega(A \cap D_1^o \dvtx m)
=
0
% \mbox{ for every $m \geq1$ and every $\omega\in\Omega$}.
\]
for every $m \geq1$ and every $\omega\in\Omega$.
However, as $A \in H_\kappa(n_\kappa(r))$
[see (\ref{hkdef})] and $\lambda(A) > 0$,
we know that $\Delta^\omega(A \cap D_1 \dvtx m) > 0$
when $\omega= \omega_\kappa(n_\kappa(r) + 1)$ and
$m = m_\kappa(n_\kappa(r) + 1)$. Thus, we arrive at a
contradiction.

Now, suppose that $(D_1^c)^o \cap A = \varnothing$. In this case,
$A \subseteq\overline{D}_1$ and with the choice of
$\omega= \omega_\kappa(n_\kappa(r) + 1)$ and
$m = m_\kappa(n_\kappa(r) + 1)$, we have
\[
\frac{\eta}{2} \lambda(A)
<
\Delta^{\omega} (A \cap D_1 \dvtx m)
=
\Delta^{\omega} (A \cap\overline{D}_1 \dvtx m)
=
\Delta^{\omega} (A \dvtx m)
\leq
\delta\lambda(A) .
\]
Here, the first inequality follows from the fact that
$A \in H_\kappa(n_\kappa(r))$ and the second follows
from (\ref{deltaa}).
Comparing the first and last
terms above, the fact that $\delta\leq\eta/12$ again yields
a contradiction.
%Choose $x_{2,1} \in C_1^o \cap A$ and $x_{2,2} \in(C_1^c)^o \cap A$.
We note that the argument above applies to any set
$A \in H_\kappa(n_\kappa(r))$ having positive Lebesgue measure.

Now, suppose that we have identified sets
$D_1,\ldots,D_l \in\mathcal{C}$, with $l \leq L-2$,
such that the join $K_l = D_1 \vee\cdots\vee D_l$
%$|K_l| = 2^l$ and $B^o \cap Q_l \neq\varnothing$ for every $B \in K_l$.
satisfies the conditions of Lemma \ref{fulljoin}.
Let $K_l = \{ B_j \dvtx j \in[2^l]\}$ and
let $x_j \in B_j^o \cap Q_l$ for each $j \in[2^l]$.
Select\vspace*{-2pt} $\alpha_{l+1} > 0$ such that for~each $j$, the interval
$I_j \stackrel{\triangle}{=}(x_j - \alpha_{l+1}, x_j + \alpha_{l+1})$
is contained in $B_j^o$ and satisfies
\[
\lambda(I_j \cap Q_l)
\geq
(1 - \varepsilon) \lambda(I_j)
=
2 \alpha_{l+1} (1 - \varepsilon) .
\]
Let $\kappa' = k_{L - l}$ and let $\{ n_{\kappa'} (l) \dvtx l \geq1 \}$
be the subsequence used to define the sub-probability $\nu_{\kappa'}$.
For each interval $I_j$,
\[
\liminf_{r \to\infty}
\lambda\bigl( I_j \cap G_{\kappa'}(n_{\kappa'} (r) ) \bigr)
\geq
\nu_{\kappa'} (I_j)
\geq
\nu_{\kappa'} (I_j \cap R_{\kappa'})
>
2 \alpha_{l+1} (1 - \varepsilon) \delta,
\]
where the last inequality follows from the previous display
and the fact that $Q_l \subseteq R_{\kappa'}$.
Choose $r$ sufficiently large so that
$\lambda( I_j \cap G_{\kappa'}(n_{\kappa'} (r)) ) >
2 \alpha_{l+1} (1 - \varepsilon) \delta$ for each $j$,
and $2^{ - n_{\kappa'}(r) } < \delta\alpha_{l+1} / 4$.

By applying the proof of Lemma \ref{Alem} to each interval $I_j$,
it is easy to see that there exist
sets $A_j \in H_{\kappa'} (n_{\kappa'}(r))$ such that
$A_j \subseteq I_j \subseteq B_j^o$, $\lambda(A_j \cap Q_{l+1}) > 0$ and
$A_j \subseteq Q_{l+1}$.
Let $D_{l+1} = C_{\kappa'}(n_{\kappa'} (r)+1) \in\mathcal{C}$. Arguments
identical to the case $l=1$ above show that for each $j$, the intersections
$A_j \cap D_{l+1}^o$ and $A_j \cap(D_{l+1}^c)^o$ are nonempty. This
completes the inductive step, and hence the proof, of Lemma \ref{fulljoin}.

Given any two dyadic intervals, they are disjoint, intersect
at one point or one contains
the other. Therefore, among the sets $D_1,\ldots,D_{L-1}$ of Lemma
\ref{fulljoin}, at most one can be a dyadic interval; the remainder
are contained in $\mathcal{C}_0$ and together have a full
join whose cells have positive Lebesgue measure. This completes
the proof of Proposition \ref{ErgVCweak}.
\end{pf*}

%As the integer $L$ above was arbitrary, Lemma \ref{join}
%implies that the VC-dimension of $\C$ is infinite, and Theorem

%s3 ###
\section[Reductions and proof of Theorem 1]{Reductions and proof of Theorem \protect\ref{ErgVC}}
\label{RPT}

As noted in the \hyperref[intro]{Introduction}, Theorem \ref{ErgVC} is
derived from
Proposition \ref{ErgVCweak} via a series of three reductions.
Two of these reductions are based on the following lemmas, whose
proofs can be found in the \hyperref[Tech]{Appendix}.
The third follows from standard
results on measure space isomorphisms.
In what follows, $A \bigtriangleup B = (A \setminus B) \cup(B
\setminus A)$
is the standard symmetric difference of two sets.
\begin{lem}
\label{RedLem1}
Let $\mathbf{X}= X_1, X_2, \ldots$ be a stationary ergodic process taking
values in $(\mathcal{X},\mathcal{S})$
and let $\mathcal{C}\subseteq\mathcal{S}$ be a countable family of
sets such
that\break $\limsup_m \Gamma_m(\mathcal{C}\dvtx\mathbf{X}) > 0$ with positive
probability.
Then $\mathcal{X}$ is necessarily uncountable and there exists a stationary
ergodic process
$\tilde{\mathbf{X}} = \tilde{X}_1, \tilde{X}_2, \ldots$\break with
values in
$(\mathcal{X},\mathcal{S})$ such that $\mathbb P(\tilde{X}_i = x) =
0$ for each $x \in\mathcal{X}$
and\break $\limsup_m \Gamma_m(\mathcal{C}\dvtx\tilde{\mathbf{X}}) > 0$ with
positive probability.
\end{lem}
\begin{lem}
\label{RedLem3}
Let $\mathcal{C}=\{ C_1, C_2,\ldots\}$ be a countable collection of Borel
subsets of $[0,1]$ such that the maximum diameter of the elements of the
join $J_n = \bigvee_{i=1}^{n} C_i$ tends to zero as $n\to\infty$.
There then exists a Borel measurable map $\phi\dvtx[0,1] \to[0,1]$
and a Borel set $V_1\subseteq[0,1]$ of measure one
such that:
\textup{(i)} $\phi$ preserves Lebesgue measure and is one-to-one on $V_1$;
\textup{(ii)} the image $V_2 = \phi(V_1)$ and the inverse map
$\phi^{-1} \dvtx V_2 \to V_1$ are Borel measurable;
\textup{(iii)} $\phi^{-1}$ preserves Lebesgue measure;
\textup{(iv)} for every set $C \in\mathcal{C}$, there is a set $U(C)$, equal
to a finite union of intervals, such that
$\lambda(\phi(C) \bigtriangleup U(C)) = 0$.
\end{lem}

%s3.1 ###
\subsection[Proof of Theorem 1]{Proof of Theorem \protect\ref{ErgVC}}

We establish the contrapositive of Theorem \ref{ErgVC} via a reduction
to Proposition \ref{ErgVCweak}.
Suppose that $\limsup_m \Gamma_m(\mathcal{C}\dvtx\mathbf{X}) > 0$ with positive
probability. Let $\mu(\cdot)$ denote the one-dimensional marginal distribution
of $\mathbf{X}$.
By Lemma \ref{RedLem1}, we may restrict our attention to the case in which
$\mu(\cdot)$ is nonatomic and $\mathcal{X}$ is uncountable.
It then follows from standard measure space isomorphism results \cite{Royd88}
that there exist Borel measurable
sets $\mathcal{X}_0 \subseteq\mathcal{X}$ and $I_0 \subseteq[0,1]$ with
$\mu(\mathcal{X}_0) = \lambda(I_0) = 1$
and an invertible map $\psi\dvtx\mathcal{X}_0 \to I_0$ such that $\psi
$ and
$\psi^{-1}$ are measurable with respect to the restricted
sigma-algebras $\mathcal{S}\cap\mathcal{X}_0$ and $\mathcal{B}\cap
I_0$, respectively, and
$\mu(A) = \lambda(\psi(A))$ for each $A \in\mathcal{S}\cap
\mathcal{X}_0$.
The event $E = \{ X_i \in\mathcal{X}_0^c \mbox{ for some } i \geq1
\}$
has probability zero, so by removing $E$ from the underlying sample space,
we may assume that $X_i(\omega) \in\mathcal{X}_0$ for each sample point
$\omega$ and each $i \geq1$.

Define $Y_i = \psi(X_i)$ for $i \geq1$ and let
$\mathcal{C}_1 = \{ \psi(C \cap\mathcal{X}_0) \dvtx C \in\mathcal{C}\}
$ be the (Borel) images in $[0,1]$
of the elements of $\mathcal{C}$.
The process $\mathbf{Y}= Y_1, Y_2, \ldots$ is stationary and ergodic
with marginal
distribution $\lambda$. If $C_1 = \psi(C \cap\mathcal{X}_0)$ is an element
of $\mathcal{C}_1$,
then $\lambda(C_1) = \mu(C \cap\mathcal{X}_0) = \mu(C)$ as $\mu
(\mathcal{X}_0) =
1$, and
$I(Y_i \in C_1) = I(\psi(X_i) \in\phi(C \cap\mathcal{X}_0)) =
I(X_i \in C)$
as $\phi(\cdot)$ is one-to-one.
Moreover, if $\mathcal{C}_1$ shatters points $u_1,\ldots,u_k \in
[0,1]$, then
$\mathcal{C}$ shatters
$\psi^{-1}(u_1), \ldots, \psi^{-1}(u_k)$.
It follows that $\Gamma_m(\mathcal{C}_1 \dvtx\mathbf{Y}) = \Gamma
_m(\mathcal{C}\dvtx\mathbf{X})$
with probability one (actually, for every $\omega$)
and that $\dim(\mathcal{C}_1) \leq\dim(\mathcal{C})$.

Let $\mathcal{C}_2 = \mathcal{C}_1 \cup\mathcal{D}$, where
$\mathcal{D}$ denotes the set of
closed dyadic subintervals of $[0,1]$. Then
$\Gamma_m(\mathbf{Y}\dvtx\mathcal{C}_2) \geq\Gamma_m(\mathbf{Y}\dvtx
\mathcal{C}_1)$ and
an easy argument shows
that \mbox{$\dim(\mathcal{D}) = 2$}. Using Lemma \ref{Sauer}
(cf. Exercise 4.1 of \cite{DevLug01}), one may show
that $\dim(\mathcal{C}_2) \leq\dim(\mathcal{C}_1) + \dim(\mathcal
{D}) + 1 \leq\dim(C_1) + 3$.
As the family $\mathcal{C}_2$ includes $\mathcal{D}$, it satisfies
the conditions of
Lemma \ref{RedLem3} above:
let $V_1,V_2$ and $\phi\dvtx[0,1] \to[0,1]$ be the associated sets and point
mapping, respectively, in the lemma. Define $Z_i = \phi(Y_i)$ for $i
\geq1$ and let
$\mathcal{C}_3 = \{ \phi(C \cap V_1) \dvtx C \in\mathcal{C}_2 \}$.
Arguments like those
above show that $\Gamma_m(\mathcal{C}_3 \dvtx\mathbf{Z}) = \Gamma
_m(\mathcal{C}_2 \dvtx\mathbf{Y})$ with
probability one
and that $\dim(\mathcal{C}_3) \leq\dim(\mathcal{C}_2)$.

By Lemma \ref{RedLem3}, for each set
$C \in\mathcal{C}_3$, there is a set $U(C)$ that is equal to a finite
union of
intervals
and is such that $\lambda(C \triangle U(C)) = 0$.
Let $\mathcal{U} = \{ U(C) \dvtx C \in C_3 \}$.
Then $\Gamma_m(\mathcal{U} \dvtx\mathbf{Z}) = \Gamma_m(\mathcal{C}_3
\dvtx\mathbf{Z})$ with probability
one and it follows from the other relations established above that
$\limsup_m \Gamma_m(\mathcal{U} \dvtx\mathbf{Z}) > 0$ with positive
probability.
Fix $L \geq1$. By Proposition \ref{ErgVCweak},
there exist sets $U(C_1),\ldots,U(C_L) \in\mathcal{U}$
such that their join has $2^L$ cells and
each cell has positive probability.
It follows that the join $J_{L} = C_1 \vee\cdots\vee C_{L}$ is also full.
As $L$ was arbitrary, Lemma \ref{join} implies that $\mathcal{C}_3$ has
infinite VC dimension,
and the same is therefore true of $\mathcal{C}$.
This completes the proof of Theorem \ref{ErgVC}.

%s4 ###
\section{Proof of VC-major and VC-graph results}
\label{PRCF}

%s4.1 ###
\subsection[Proof of Proposition 1]{Proof of Proposition \protect\ref{VCMajor}}
\label{VCMajorProof}

Let $\mathbf{X}$ be a stationary ergodic process.
Suppose, first, that $\mathcal{F}$ is bounded, with constant envelope
$M <
\infty$.
Fix $\varepsilon> 0$ and select an integer $K$ such that
$2 M / K \leq\varepsilon$. For each $f \in\mathcal{F}$, define the
approximation
\[
\overline{f}(x)
=
M - \frac{2M}{K} \sum_{j=1}^K I\bigl( f(x) \leq M - 2M j / K\bigr) .
\]
Note that $\overline{f}(x) - \varepsilon\leq f(x) \leq\overline{f}(x)$
for each
$x \in\mathcal{X}$ and thus, by an elementary bound,
\[
0
\leq
\Gamma_m(\mathcal{F}\dvtx\mathbf{X})
\leq
2 \varepsilon
+ \Gamma_m(\overline{\mathcal{F}}\dvtx\mathbf{X}) ,
\]
where $\overline{\mathcal{F}} = \{ \overline{f} \dvtx f \in\mathcal
{F}\}$. It follows readily
from Theorem \ref{ErgVC} and the assumption that $\dim_{\mathrm{VC}}(\mathcal{F})$
is finite that $\Gamma_m(\overline{\mathcal{F}}\dvtx\mathbf{X}) \to0$
with probability
one as $n$ tends to infinity. As $\varepsilon> 0$ was arbitrary,
we conclude that $\Gamma_m(\mathcal{F}\dvtx\mathbf{X}) \to0$ with
probability one
as well.

Now, suppose that $\mathcal{F}$ has an envelope $F$ such that $EF(X) <
\infty
$. Fix
$0 < M < \infty$ and for each $f \in\mathcal{F}$, define $f_M(x) =
f(x) I(
F(x) \leq M)$.
Let $\mathcal{F}_M = \{ f_M \dvtx f \in\mathcal{F}\}$. Then, by an
elementary bound and an
application
of the ergodic theorem to $F(x) I( F(x) \leq M)$,
we have
\[
0
\leq
\limsup_{m \to\infty} \Gamma_m(\mathcal{F}\dvtx\mathbf{X})
\leq
\limsup_{m \to\infty} \Gamma_m(\mathcal{F}_M \dvtx\mathbf{X})
+
2 E\bigl[F(X)I\bigl(F(X) >M\bigr)\bigr] .
\]
A straightforward argument shows that $\mathcal{F}_M$ is a VC-major
class and therefore,
by the result above,
the first term on the right-hand side is equal to zero. The second term
can be
made arbitrarily small by choosing $M$ sufficiently large.

%s4.2 ###
\subsection[Proof of Proposition 2]{Proof of Proposition \protect\ref{VCGraph}}
\label{VCGraphProof}

Let $\mathbf{X}$ be a stationary ergodic process with one-dimensional
marginal distribution $\mu$.
Let $M < \infty$ be an envelope for $\mathcal{F}$. Replacing each $f
\in\mathcal{F}$
by $(f + M)/2M$, we may assume
without loss of generality that each $f \in\mathcal{F}$ takes values in
$[0,1]$ and,
therefore,
\[
G_f = \{ (x,s) \dvtx x \in\mathcal{X}\mbox{ and } 0 \leq s \leq f(x)
\leq1 \} .
\]
Let $Y_1, Y_2, \ldots\in[0,1]$ be independent, uniformly distributed
random variables defined on the same probability space as,
and independent of, the process $\mathbf{X}$. For $i \geq1$, define
$Z_i = (X_i,Y_i) \in\mathcal{X}\times[0,1]$.
It follows from standard results in ergodic theory (cf. \cite{Pet83}) that
the process $\mathbf{Z}= Z_1, Z_2, \ldots$ is stationary and ergodic.
Let $Z = (X,Y)$ be distributed as $Z_1$.
By an application of Fubini's theorem, for each $f \in\mathcal{F}$,
we have
%
%e18 ###
%
\begin{eqnarray}
\label{vcg1}
\mathbb P(Z \in G_f)
&=&
(\mu\otimes\lambda) (G_f)
=
\int_\mathcal{X} \lambda( (G_f)_x ) \,d\mu(x)\nonumber\\[-8pt]\\[-8pt]
&=&
\int_\mathcal{X} f(x) \,d\mu(x)
=
Ef(X),\nonumber
\end{eqnarray}
where $G_x = \{ s \dvtx(x,s) \in G \}$ denotes the $x$-section of $G$.
Moreover,
%
%e19 ###
%
\begin{equation}
\label{vcg2}
\frac{1}{m} \sum_{i =1}^m I( Z_i \in G_f )
=
\frac{1}{m} \sum_{i =1}^m I\bigl( Y_i \leq f(X_i) \bigr) .
\end{equation}
By an elementary bound,
$\Gamma_m(\mathcal{F}\dvtx\mathbf{X}) \leq\Gamma^1_m(\mathcal
{F}\dvtx\mathbf{Z}) + \Gamma^2_m(\mathcal{F}\dvtx\mathbf{Z})$,
where
\[
\Gamma^1_m(\mathcal{F}\dvtx\mathbf{Z})
=
\sup_{f \in\mathcal{F}}
\Biggl| \frac{1}{m} \sum_{i =1}^m I\bigl(Y_i \leq f(X_i)\bigr) - Ef(X) \Biggr|
\]
and
\[
\Gamma^2_m(\mathcal{F}\dvtx\mathbf{Z})
=
\sup_{f \in\mathcal{F}}
\Biggl| \frac{1}{m} \sum_{i =1}^m
\bigl[ I\bigl(Y_i \leq f(X_i)\bigr) - f(X_i) \bigr] \Biggr| .
\]
It follows from (\ref{vcg1}) and (\ref{vcg2}) that
\[
\Gamma^1_m(\mathcal{F}\dvtx\mathbf{Z})
=
\sup_{G \in\mathcal{G}}
\Biggl| \frac{1}{m} \sum_{i =1}^m I(Z_i \in G)
- \mathbb P(Z \in G) \Biggr| ,
\]
which tends to zero with probability one by Theorem \ref{ErgVC}
and the assumption that $\mathcal{G}$ is a VC class.
To analyze the second supremum, note that when
$X_1= x_1,\ldots, X_m = x_m$ are fixed,
\[
\Gamma^2_m(\mathcal{F}\dvtx(x_1,Y_1),\ldots,(x_n,Y_n))
=
\sup_{f \in\mathcal{F}}
\Biggl| \frac{1}{m} \sum_{i =1}^m \bigl[ I\bigl(Y_i \leq f(x_i)\bigr) - \mathbb P\bigl(Y_i \leq
f(x_i)\bigr) \bigr] \Biggr|
\]
and that $Y_1,\ldots,Y_n$ remain independent under
this conditioning. By a routine modification of standard
empirical process arguments like those in Theorem 3.1
of Devroye and Lugosi \cite{DevLug01}, one may establish that
\[
E[ \Gamma^2_m(\mathcal{F},\mathbf{Z}) | X_1^n ]
\leq
2 \biggl( \frac{\ln2 S_m(\mathcal{G})}{m} \biggr)^{1/2}
\stackrel{\triangle}{=}
L_m .
\]
Here, $S_m(\mathcal{G})$ is the (maximal) shatter coefficient of
$\mathcal{G}$ defined by
\[
S_m(\mathcal{G})
=
\max
\bigl|\bigl\{ G \cap\{z_1,\ldots,z_m\} \dvtx G \in\mathcal{G}\bigr\}\bigr|,
\]
where the maximum is taken over all $m$-sequences
$z_1,\ldots,z_m \in\mathcal{X} \times[0,1]$. As $\mathcal{G}$ has finite
VC dimension, $V$ say, it follows from Sauer's Lemma \ref{Sauer}
above that $S_m(\mathcal{G}) \leq(m+1)^V$
and, consequently, that $L_m = O( (\ln m / m)^{1/2} )$. A straightforward
application of McDiarmid's bounded difference inequality
(cf. Theorem 2.2 of \cite{DevLug01}) shows that for $t > 0$,
\[
\mathbb P\bigl( \Gamma^2_m(\mathcal{F}\dvtx\mathbf{Z}) \geq L_m + t | X_1^n \bigr)
\leq
e^{-2 m t^2} .
\]
Taking expectations, the same bound holds for the unconditional
probability and it then
follows from a simple application of the first Borel--Cantelli lemma that
$\Gamma^2_m(\mathcal{F}\dvtx\mathbf{Z})$
tends to zero with probability one as $m$ tends to infinity.

\begin{appendix}\label{Tech}
%s5 ###
\section*{Appendix}

%s5.1 ###
\subsection{\texorpdfstring{Proof of Lemma \protect\ref{RedLem1}}{Proof of Lemma 5}}
\label{R1}

Following arguments like those in Breiman \cite{Breim92},
we may assume, without loss of generality, that $\mathbf{X}= \{ X_i
\dvtx
-\infty
< i < \infty\}$
is a two-sided process and that $\mathbf{X}$ is defined on a
probability space
$(\Omega, \mathcal{F}, \mathbb P)$ via a left shift transformation
and a projection
map. Specifically,
$\Omega$
%= \times_{i = -\infty}^\infty\X$
is the set of all bi-infinite
sequences $\omega= (\omega_i)_{i = -\infty}^\infty$, where
$\omega_i \in\mathcal{X}$ for each $i$, and $\mathcal{F}= \bigotimes
_{i = -\infty
}^\infty\mathcal{S}$ is the usual
product sigma-field.
%on $\Omega$, generated by the finite dimensional cylinder sets
%$E = \{ \omega\in\Omega\dvtx\omega_l \in A_l \mbox{ for } |l| \leq k
%with $A_l \in\sigS$ and $k \geq1$.
We may further assume that
$X_i(\omega) = X_0(T^i \omega)$, where $X_0\dvtx\Omega\to\mathcal{X}$
is the coordinate
projection $X_0(\omega) = \omega_0$ and $T\dvtx\Omega\to\Omega$ is
the standard
left-shift transformation defined by $(T \omega)_i = \omega_{i-1}$.
The stationarity of $\mathbf{X}$ implies that
%$\mathbb P(E) = \mathbb P(TE) = \mathbb P(T^{-1}E)$ for each $E \in
$T$ and $T^{-1}$ preserve $\mathbb P(\cdot)$. Ergodicity
of $\mathbf{X}$ ensures that $T$ is ergodic: if $TA = A$, then
$\mathbb P(A) = 0$
or $1$.

As noted by Steele \cite{Stee78}, the subadditive ergodic theorem implies
that the random variables $\Gamma_m(\mathcal{C}\dvtx\mathbf{X})$
converge with
probability one to a constant.
In particular, if $\limsup_m \Gamma_m(\mathcal{C}\dvtx\mathbf{X}) > 0$
with positive
probability, then it follows that
%
%e20 ###
%
\begin{equation}
\label{set}
\liminf_m \Gamma_m(\mathcal{C}\dvtx\mathbf{X}) > 0
\qquad\mbox{with probability one} .
\end{equation}
This stronger converse of the Glivenko--Cantelli property will be
needed in what follows.

Let $A = \{x \in\mathcal{X}\dvtx\mu(\{x\}) = 0 \}$ contain the
nonatomic points
of $\mathcal{X}$.
If $A^c = \varnothing$, then $\mathcal{X}$ is uncountable and there is nothing
else to
prove. Assume, then, that $A^c \neq\varnothing$.
As $A^c$ consists of the (finite or countable) set of points in
$\mathcal{X}$
having positive
$\mu$-measure, it follows that $A \in\mathcal{S}$. Given $\varepsilon> 0$,
we may express
$A^c$ as a disjoint union $A_1 \cup A_2$ such that the cardinality of $A_1$
is finite and $\mu(A_2) < \varepsilon$.
Let $\hat{\mu}_m(A) = m^{-1} \sum_{i=1}^m I(X_i \in A)$ denote the empirical
measure of $X_1,\ldots,X_m$. By an elementary bound,
\begin{eqnarray*}
\Gamma_m(\mathcal{C}\dvtx\mathbf{X})
& \leq&
\Gamma_m(\mathcal{C}\cap A \dvtx\mathbf{X})
% | \hat{\mu}_m(C \cap A) - \mu(C \cap A) |
+
{\sum_{x \in A_1}}
| \hat{\mu}_m(\{x\}) - \mu(\{x\}) |
+
\hat{\mu}_m(A_2) + \mu(A_2) .
\end{eqnarray*}
As $m$ tends to infinity, the second term above tends to zero and the last
two terms are together less than $2 \varepsilon$. As $\varepsilon> 0$ was
arbitrary, we conclude that $\mu(A) > 0$ so that $\mathcal{X}$ is uncountable.
Moreover, (\ref{set}) implies that
$\liminf_m \Gamma_m(\mathcal{C}\cap A \dvtx\mathbf{X}) > 0$ with
probability one.

Let $\Omega_A$ denote the set of $\omega\in\Omega$
such that $\omega_0 \in A$ and both index sets
$\{ i \geq1 \dvtx w_i \in A\}$ and $\{ i \leq-1 \dvtx w_i \in A\}$ are infinite.
By the ergodic theorem, $\mathbb P(\Omega_A) = \mu(A) > 0$.
For $\omega\in\Omega_A$, define
$\tau(\omega) = \min\{ k \geq1\dvtx T^k \omega\in A \}$ (which is finite)
and the induced
transformation $\tilde{T}\dvtx\Omega_A \to\Omega_A$ by
$\tilde{T} \omega= T^{ \tau(\omega) } \omega$.
Routine arguments from ergodic theory \cite{Pet83} show that
$\tilde{T}$ is invertible, is measurable on the restricted
sigma-field $\mathcal{F}_A = \mathcal{F}\cap\Omega_A$, preserves the
normalized measure $\mathbb P_A (\cdot) = \mathbb P(\cdot) / \mathbb
P(\Omega_A)$ on
$(\Omega_A,\mathcal{F}_A)$ and is ergodic.
For the sake of completeness, we provide a sketch of the proofs using
a geometric argument from ergodic theory
known as the \textit{Kakutani skyscraper}. For each positive
integer $k$, define $A_k = \{ \omega\in\Omega_A\dvtx\tau(\omega)=k \}$.
The sets $A_1, A_2,\ldots$ then partition $\Omega_A$.
Moreover, $\bigcup_{k=1}^{\infty} \bigcup_{i=0}^{k-1} T^i A_k$ is a disjoint
union containing almost every point in $\Omega$.
The Kakutani skyscraper of $\Omega_A$ is created by
stacking the sets $T^1 A_k,\ldots, T^{k-1} A_k$ above
$A_k$ for each $k \geq1$.

The measurability
of $\tilde{T}$ follows from the fact that each $A_k$ is measurable
and that $\tilde{T}$ restricted to $A_k$ equals $T^k$ restricted to $A_k$.
Invertibility of $\tilde{T}$ follows directly from the invertibility
of $T$ and the construction of the\vspace*{1pt} Kakutani skyscraper.
In particular, let $\omega_1 \neq\omega_2$ be points in $\Omega_A$.
Then $\tilde{T}(\omega_1) = T(T^{\tau(\omega_1)-1}\omega_1)$ and
$\tilde{T}(\omega_2) = T(T^{\tau(\omega_2)-1}\omega_2)$.
As $T$ is invertible, and $T^{\tau(\omega_1)-1}(\omega_1)$ and
$T^{\tau(\omega_2)-1}(\omega_2)$ are distinct points in the
Kakutani skyscraper,
it follows that
$\tilde{T}(\omega_1) \neq\tilde{T}(\omega_2)$.
The measure-preserving property of $\tilde{T}$ follows from the fact that
$T$ is measure preserving on each of the sets $A_k$.
To establish ergodicity, suppose that
$B \subset\Omega_A$ is a set of positive measure that is invariant
for $\tilde{T}$.
The set $C = \bigcup_{i=-\infty}^{\infty} T^{i}B$ is
invariant under $T$, and $C \cap A = B$
since $B$ is invariant for $\tilde{T}$.
As $T$ is ergodic, $C$ contains $A$ and so
$A=B$. It follows that $\tilde{T}$ is ergodic.

Define $\tilde{X}_0\dvtx\Omega_A \to\mathcal{X}$ by $\tilde
{X}_0(\omega) =
\omega_0$
and $\tilde{X}_i(\omega) = \tilde{X}_0(\tilde{T}^i \omega)$ for
$-\infty< i < \infty$.
The process $\tilde{\mathbf{X}} = \{ \tilde{X}_i \}$
defined on $(\Omega_A, \mathcal{F}_A, \mathbb P_A)$ is
then stationary and ergodic, takes values in $(\mathcal{X}, \mathcal
{S})$ and has
marginal distribution $\mu_A(\cdot) = \mu(\cdot) / \mu(A)$ with no
point masses.

We wish to show that
$\limsup_m \Gamma_m(\mathcal{C}\dvtx\tilde{\mathbf{X}}) > 0$ with
positive $\mathbb P_A$-probability.
To this end, for each $\omega\in\Omega_A$, define $\tau_0(\omega)
= 0$,
$\tau_1(\omega) = \tau(\omega)$ and
$\tau_{l+1}(\omega) = \min\{ k > \tau_{l}(\omega) \dvtx\omega_k \in
A \}$.
By definition of $\Omega_A$, each function $\tau_l$ is finite.
For each $m \geq1$, $C \in\mathcal{C}$ and $\omega\in\Omega_A$,
%
%e21 ###
%
\begin{eqnarray}\label{xtildex}\qquad
\frac{1}{m} \sum_{i=0}^{m-1} I\bigl( \tilde{X}_i (\omega) \in C \bigr)
& = &
\frac{1}{m}
\sum_{j=0}^{\tau_{m-1}(\omega)} I\bigl( X_j (\omega) \in C \cap A \bigr)
\nonumber\\[-8pt]\\[-8pt]
& = &
\frac{1}{\mu(A)} W_m(\omega)
\frac{1}{\tau_{m-1}(\omega)}
\sum_{j=0}^{\tau_{m-1}(\omega)}
I\bigl( X_j (\omega) \in C \cap A \bigr),\nonumber
\end{eqnarray}
where we have defined $W_m = \mu(A) \tau_{m-1} / m$.
By the ergodic theorem, for $\mathbb P_A$-almost every $\omega\in
\Omega_A$,
\[
\frac{m}{\tau_{m-1}(\omega)}
=
\frac{1}{\tau_{m-1}(\omega)}
\sum_{j=0}^{\tau_{m-1}(\omega)} I\bigl( X_j (\omega) \in C \cap A \bigr)
\to\mu(A)
\]
as $m$ tends to infinity
and so $W_m \to1$ with $\mathbb P_A$
probability one. Omitting the dependence on $\omega$,
it follows from (\ref{xtildex}) and the definition of $\mu_A(\cdot)$
that
\begin{eqnarray*}
\Gamma_m(\mathcal{C}\dvtx\tilde{\mathbf{X}})
& = &
\sup_{C \in\mathcal{C}}
\Biggl|
\frac{1}{m} \sum_{i=0}^{m-1} I( \tilde{X}_i \in C )
- \mu_A(C)
\Biggr| \\
& = &
\frac{1}{\mu(A)}
\sup_{C \in\mathcal{C}}
\Biggl|
W_m
\frac{1}{\tau_{m-1}}
\sum_{j=0}^{\tau_{m-1}}
I( X_j \in C \cap A )
- \mu(C \cap A)
\Biggr| \\
& \geq&
\frac{1}{\mu(A)}
\Gamma_{\tau_{m-1}} (\mathcal{C}\cap A \dvtx\mathbf{X})
-
|W_m-1|
\sup_{C \in\mathcal{C}}
\Biggl|
\frac{1}{\tau_{m-1}}
\sum_{j=0}^{\tau_{m-1}}
I( X_j \in C \cap A )
\Biggr| \\
& \geq&
\frac{1}{\mu(A)}
\Gamma_{\tau_{m-1}} (\mathcal{C}\cap A \dvtx\mathbf{X})
-|W_m-1| .
\end{eqnarray*}
The first inequality above follows by writing
$W_m = 1 + (W_m - 1)$ and then using the
elementary bound
${\sup_\alpha}|a_\alpha- b_\alpha|
\geq{\sup_\alpha}|a_\alpha| - {\sup_\alpha}|b_\alpha|$.
It follows from the last display that
\begin{eqnarray*}
\limsup_m \Gamma_m(\mathcal{C}\dvtx\tilde{\mathbf{X}})
& \geq&
\liminf_m \Gamma_m(\mathcal{C}\dvtx\tilde{\mathbf{X}}) \\
& \geq&
\frac{1}{\mu(A)}
\liminf_m \Gamma_{\tau_{m-1}}(\mathcal{C}\cap A \dvtx\mathbf{X}) \\
& \geq&
\frac{1}{\mu(A)}
\liminf_m \Gamma_{m}(\mathcal{C}\cap A \dvtx\mathbf{X})
\end{eqnarray*}
and the argument above shows that the final term
is positive with $\mathbb P_A$-probability one. This
completes the proof.

%s5.2 ###
\subsection{\texorpdfstring{Proof of Lemma \protect\ref{RedLem3}}{Proof of Lemma 6}}
\label{R3}

The isomorphism $\phi$ is defined as a limit of isomorphisms $\phi_n$.
The maps $\phi_n$ are defined inductively. To begin, let
\begin{eqnarray*}
\phi_1(x) =
\cases{
\lambda([0,x] \cap C_1) , &\quad if $x \in C_1 $, \cr
\lambda(C_1)+\lambda([0,x] \cap C_1^c) , &\quad if $x \in C_1^c $.}
\end{eqnarray*}
Then $\phi_1$ maps $C_1$ into $[0,\lambda(C_1)]$ and $C_1^c$ into
$[\lambda(C_1),1]$.
By standard arguments, $\phi_1$ is Lebesgue measure-preserving
and a bijection almost everywhere.

Suppose now that maps $\phi_1,\ldots,\phi_n$ have been defined in
such a way that: (i)~for each element $A$ of the
join $J_n = \bigvee_{i=1}^{n} C_i$ and each
$x \in A$, $\phi_n(x) = \beta_n(A) + \lambda([0,x] \cap A)$,
where $\beta_n(A)$ is a constant; (ii) the intervals
$\{ [\beta_n(A),\break\beta_n(A)+\lambda(A)) \dvtx A \in J_n \}$ form a disjoint
covering of $[0,1)$.
For each each $A \in J_n$ and each $x \in A$,
define
\begin{eqnarray*}
\phi_{n+1}(x) =
\cases{
\beta_n(A) + \lambda( [0,x] \cap A \cap C_{n+1}) ,
&\quad if $x \in A \cap C_{n+1} $, \cr
\beta_n(A) + \lambda( A \cap C_{n+1}) + \lambda( [0,x] \cap A \cap
C_{n+1}^c) ,
&\quad if $x \in A \cap C_{n+1}^c $.}
\end{eqnarray*}
With these definitions, properties (i) and (ii) hold for $J_{n+1}$
and $\phi_{n+1}$.
%$\{ [\beta_{n+1}(A),\beta_{n+1}(A)+\lambda(A)) \dvtx A\in J_{n+1} \}$
%is a disjoint covering of $[0,1)$ and
Moreover, $\phi_1, \phi_2,\ldots$ have the property that
for each $n$, each cell $A \in J_n$ and each $m \geq n$,
the function $\phi_m$ is a Lebesgue-measure-preserving almost
everywhere bijection from $A$
into $[\beta_n(A),\beta_n(A) + \lambda(A)]$.
In particular,
for each $A \in J_n$ and each $m \geq n$,
\[
\operatorname{cl}(\phi_n(A))
=
\operatorname{cl}(\phi_m(A))
=
[\beta_n(A),\beta_n(A) + \lambda(A)] ,
\]
where $\operatorname{cl}(U)$ denotes the closure of $U$.

Fix $x \in[0,1]$ for the moment and, for $n \geq1$,
let $A_n(x)$ be the cell of $J_n$ containing $x$.
Note that the sequence $\phi_n(x), \phi_{n+1}(x), \ldots$ is
contained in the
interval $\operatorname{cl}(\phi_n(A_n(x)))$, whose diameter is equal to
$\lambda(A_n(x)) \leq\operatorname{diam}(A_n(x))$. By assumption, the
latter quantity tends to zero as $n \to\infty$ and so
$\{ \phi_n(x) \dvtx n \geq1 \}$ is a Cauchy sequence. Let $\phi(x)$
denote its limit. Then $\phi(\cdot)$ is a limit of measurable
functions, hence measurable.

We claim that $\operatorname{cl}(\phi(A)) = \operatorname{cl}(\phi
_n(A))$ for every
$n \geq1$ and every $A \in J_n$.
To see this, fix $A \in J_n$.
If $y \in\operatorname{cl}(\phi(A))$, then there exist
$x_1,x_2,\ldots\in A$
such that $\phi(x_m) \to y$. By definition of $\phi(\cdot)$, there exist
integers $r_1, r_2,\ldots$ tending to infinity
such that $\phi_{r_m}(x_m) \to y$. As each value
$\phi_{r_m}(x_m) \in\phi_{r_m}(A) \subseteq\operatorname{cl}(\phi_n(A))$,
we have $y \in\operatorname{cl}(\phi_n(A))$. Thus,
$\operatorname{cl}(\phi(A)) \subseteq\operatorname{cl}(\phi
_n(A))$, the latter set being equal
to the interval $I_A = [\beta_n(A),\beta_n(A) + \lambda(A)]$.
To establish the converse,
let $y_0 \in I_A^o$ and $\varepsilon> 0$ be such that
$(y_0 - \varepsilon, y_0 + \varepsilon) \subseteq I_A$.
By the shrinking diameter assumption on the
joins $J_m$, there exists an integer $m$ and a cell
$A' \in J_m$ such that
$A' \subseteq A$ and $\operatorname{cl}(\phi_m(A')) \subseteq I_0$
has positive
measure. Thus, if $x \in A'$, then $\phi_r(x) \in\operatorname
{cl}(\phi_m(A'))$
for $r \geq m$
and, therefore, $\phi(x) \in I_0$.
As $\varepsilon> 0$ was arbitrary, it follows that
$I_A^o \subseteq\operatorname{cl}(\phi(A))$ and, consequently,
$I_A \subseteq\operatorname{cl}(\phi(A))$ as well.

We now establish that the map $\phi$ preserves Lebesgue measure.
To this end, for each $n \geq1$, we define
\[
\mathcal{Q}_n
=
\{ \operatorname{cl}(\phi(A)) \dvtx A \in J_n \}
\cup
\bigl\{ \{ \beta_n(A) \} \dvtx A \in J_n \bigr\} \cup\{ \{1\} \}
\]
to be the collection of intervals into which the elements of $J_n$
are mapped and the endpoints of these intervals. We wish to show that
$\lambda(\phi^{-1} B) = \lambda(B)$ for each \mbox{$B \in\mathcal{Q}_n$}.
First, suppose that $\alpha$ is the endpoint of some
interval $\operatorname{cl}(\phi(A'))$ with \mbox{$A' \in J_n$}. Fix
$\varepsilon> 0$ and let
$m \geq n$ be
large enough so that $\max\{ \lambda(A) \dvtx A \in J_m \} \leq
\varepsilon/2$.
Let $A_1,\ldots,A_r$ be those elements of $J_m$ such that
$\operatorname{cl}(\phi_m(A_j))$ contains the point $\alpha$. Then
$\phi^{-1} \{\alpha\} \subseteq\bigcup_{j=1}^r A_j$ and at most
two of the sets $A_j$ can have positive measure. It follows
that $\lambda(\phi^{-1} \{\alpha\}) \leq\varepsilon$, and, as
$\varepsilon> 0$ was arbitrary, we have $\lambda(\phi^{-1} \{\alpha\})
= 0$.
Now, suppose that $B \in\mathcal{Q}_n$ is of the form
$B = \operatorname{cl}(\phi(A)) = [\alpha_1,\alpha_2]$
for some element $A \in J_n$.
Then $\phi^{-1} B = A \cup\phi^{-1} \{ \alpha_1 \} \cup\phi^{-1}
\{ \alpha_2 \}$
and therefore
\[
\lambda(\phi^{-1} B)
=
\lambda(A)
=
\lambda(\operatorname{cl}(\phi_n(A)))
=
\lambda(\operatorname{cl}(\phi(A)))
=
\lambda(B).
\]
It follows from these arguments
that $\lambda(\phi^{-1} B) = \lambda(B)$ for each $B \in\bigcup_{m
\geq1} \mathcal{Q}_m$.
As the latter collection generates the
Borel sigma-field of $[0,1]$ and is closed under intersections, $\phi$
preserves
Lebesgue measure.

Next, we show that $\phi$ is one-to-one on a Borel subset of $[0,1]$ with
full measure.
%Recall, that each element in the join $J_n = \bigvee_{i=1}^{n}C_i$ maps
%almost everywhere to an interval under the mapping $\phi(\cdot)$.
Let $\mathcal{Q}^0 = \bigcup_{m=1}^{\infty} \{ \beta_m(A) \dvtx A \in
J_m \} \cup\{ \{1\} \}$
be the (countable) set of endpoints of the
intervals $\{\operatorname{cl}( \phi(A) ) \dvtx A \in J_m, m \geq1 \}$.
Since $\phi^{-1}$ preserves Lebesgue measure,
$\lambda(\phi^{-1}\mathcal{Q}^0) = 0$. Define
$V_1 = [0,1] \setminus\phi^{-1}\mathcal{Q}^0$, so that $\lambda(V_1)=1$.
Let $x_1$ and $x_2$ be distinct points
in $V_1$. Since the diameters of the elements of
$J_n$ tend to zero, there exists an $n$ such that
$x_1$ and $x_2$ are contained in different elements of $J_n$.
Thus, $\phi_n$ maps $x_1$ and
$x_2$ to distinct intervals, which may intersect only at their endpoints.
Hence, $\phi$ also maps $x_1$ and
$x_2$ to distinct intervals. Since $V_1$ excludes points that map
to endpoints of these intervals,
$\phi(x_1) \neq\phi(x_2)$.
Therefore, $\phi$ is a bijection on $V_1$ and we have established
conclusion (i) of the lemma.

Conclusion (ii) of the lemma follows from (i) and general results
concerning measurable maps of complete separable metric
spaces; see Corollary 3.3 of Parthasarathy \cite{Par67}.
To establish (iii), note that for any measurable subset $A \subseteq V_1$,
$\lambda(\phi(A)) = \lambda(\phi^{-1}(\phi(A))) = \lambda(A)$
since $\phi$ is measure-preserving and one-to-one on $V_1$.

To establish conclusion (iv), let $C \in\mathcal{C}$. There then exist
positive integers
$k$ and $n$ such that $C = \bigcup_{i=1}^{k} A_i$, where
$A_1, A_2, \ldots,A_k$ are (disjoint) cells in $J_n$.
Let $U(C) = \bigcup_{i=1}^{k} [\beta_n(A_i),\beta_n(A_i)+\lambda
(A_i)]$. Then
$\phi(C)=\bigcup_{i=1}^{k}\phi(A_i) \subseteq\bigcup_{i=1}^{k}$
cl($\phi(A_i)$)$=U(C)$
and $\lambda(\phi(C)) = \sum_{i=1}^{k} \lambda(A_i) = \lambda(U(C))$.
Thus, $\lambda(\phi(C)\bigtriangleup U(C)) = \lambda(U(C)\setminus
\phi(C)) = 0$.
\begin{Remark*} The condition that the cells of the joins have
diminishing diameters,
rather than measures tending to zero, is necessary.
If, for example, $C_n = \bigcup_{i=0}^{2^n-1} [\frac
{2i}{2^{n+1}},\frac{2i+1}{2^{n+1}})$
for positive integers $n$, then the limiting map is
$\phi(x) = 2x \operatorname{mod} 1$.
\end{Remark*}
\end{appendix}

% imsref loaded by lrinkeviciute, 2010-02-12 15:31:01
%

%
\printaddresses

\end{document}